\documentclass[12pt]{amsart}

\usepackage[mathscr]{eucal}
\usepackage{amsmath,amssymb,amscd,amsthm}
\usepackage{color}
\usepackage{amsmath,amsthm,amssymb,graphicx}
\usepackage{epsfig}

\theoremstyle{definition}

\theoremstyle{plane}

\def \beq{ \begin{equation} }
\def \eeq{\end{equation}}

\usepackage[top=4cm, bottom=4cm, left=3.5cm, right=3.5cm]{geometry}

\title{The classical $N$-body problem in the context of curved space}
\begin{document}
\maketitle
\markboth{Florin Diacu}{The classical $N$-body problem in the context of curved space}
\author{\begin{center}
{\bf Florin Diacu}\\
\smallskip
{\footnotesize Pacific Institute for the Mathematical Sciences\\
and\\
Department of Mathematics and Statistics\\
University of Victoria\\
P.O.~Box 3060 STN CSC\\
Victoria, BC, Canada, V8W 3R4\\
diacu@uvic.ca\\
}\end{center}

}


\begin{center}
\today
\end{center}


\begin{abstract}
We provide the differential equations that generalize the Newtonian $N$-body problem of celestial mechanics to spaces of constant Gaussian curvature, $\kappa$, for all $\kappa\in\mathbb R$. In previous studies, the equations of motion made sense only for $\kappa\ne 0$. The system derived here does more than just include the Euclidean case in the limit $\kappa\to 0$: it recovers the classical equations for $\kappa=0$. This new expression of the laws of motion allows the study of the $N$-body problem in the context of constant curvature spaces and thus offers a natural generalization of the Newtonian equations that includes the classical case. We end the paper with remarks about the bifurcations of the first integrals.
\end{abstract}

\vspace{-0.1cm}

\section{Introduction}

The idea that geometry and physics are intimately related made its way in human thought during the early part of the 19th century. Gauss measured the angles of a triangle formed by three mountain peaks near G\"ottingen, apparently hoping to learn whether the universe has positive or negative curvature, but the inevitable observational errors rendered his results inconclusive, \cite{Diacu03}, \cite{Goe},  \cite{Halsted}. In the 1830s, Bolyai and Lobachevsky took these investigations further. They independently addressed the connection between geometry and physics by seeking a natural extension of gravitation from Euclidean to hyperbolic space, \cite{Bol}, \cite{Lob}. Their idea led to the study of the Kepler problem and the 2-body problem in spaces of nonzero constant Gaussian curvature, $\kappa\ne 0$, two problems that are not equivalent, unlike in Euclidean space, \cite{Shc}. A detailed history of the results obtained in this direction since Bolyai and Lobachevsky, as well as the reasons why their approach provides a natural way of extending gravitation to spaces of constant Gaussian curvature (an aspect we also briefly mention in Section 2), can be found in \cite{Diacu03}, \cite{Diacu05}, and \cite{Diacu07}.

In some recent studies, such as \cite{Diacu00}, \cite{Diacu01}, \cite{Diacu02}, \cite{Diacu03}, \cite{Diacu04}, \cite{Diacu05}, \cite{Diacu06}, \cite{Diacu07}, \cite{Diacu08}, \cite{Diacu09}, \cite{Diacu10}, \cite{Diacu11}, \cite{Diacu12}, \cite{Diacu13}, \cite{Diacu14}, \cite{Perez}, we introduced a suitable framework for generalizing the equations of motion suggested by Bolyai and Lobachevsky to $N\ge 2$ bodies. Like the curved Kepler problem and the curved 2-body problem, our equations made sense in spaces of constant Gaussian curvature $\kappa\ne 0$,  i.e., on 3-spheres of radius $R=\kappa^{-1/2}$ embedded in $\mathbb R^4$, for $\kappa>0$, and on hyperbolic 3-spheres of imaginary radius $iR=\kappa^{-1/2}$ embedded in the Minkowski space $\mathbb R^{3,1}$, for $\kappa<0$. 
But whether written in extrinsic or intrinsic coordinates, 
these equations contain undetermined expressions for $\kappa=0$, although we can recover the classical Newtonian system when $\kappa\to 0$. So a study of the flat case in the context of curved space, including some understanding of the bifurcations and the stability of solutions when the parameter $\kappa$ is varied through 0, was impossible to perform in that setting.

In this paper we derive some equations of motion that overcome the difficulties mentioned above. Using a coordinate system in $\mathbb R^4$ having the origin at the North-Pole of the 3-spheres, we prove that the $N$-body problem in spaces of constant Gaussian curvature $\kappa\in\mathbb R$ can be written as 
\begin{equation}\label{new}
\ddot{\bf r}_i=\sum_{j=1, j\ne i}^N\frac{m_j\Big[{\bf r}_j-\Big(1-\frac{\kappa r_{ij}^2}{2}\Big){\bf r}_i+\frac{r_{ij}^2{\bf r}}{2}\Big]}{r_{ij}^3\Big(1-\frac{\kappa r_{ij}^2}{4}\Big)^{3/2}}-(\dot{\bf r}_i\cdot\dot{\bf r}_i)(\kappa{\bf r}_i+{\bf r}), \ \  i=\overline{1,N}, 
\end{equation}
where $m_1,m_2,\dots,m_N>0$ represent the masses, $\sigma=+1$ for $\kappa\ge 0$ and $\sigma=-1$ for $\kappa<0$,
$$
{\bf r}=(0,0,0,\sigma|\kappa|^{1/2}), \ \ {\bf r}_i=(x_i, y_i, z_i,\omega_i), \ i=\overline{1,N},
$$ 
the dot $\cdot$ denotes the standard inner product of signature $(+,+,+,+)$ for $\kappa\ge 0$, but the Lorentz inner product of signature $(+,+,+,-)$ for $\kappa<0$, and
$$
r_{ij}:=
\begin{cases}
[(x_i-x_j)^2+(y_i-y_j)^2+(z_i-z_j)^2+(\omega_i-\omega_j)^2]^{1/2}\ \
{\rm for}\ \ \kappa>0\cr
[(x_i-x_j)^2+(y_i-y_j)^2+(z_i-z_j)^2]^{1/2}\ \hspace{2.45cm}
{\rm for}\ \ \kappa=0\cr
[(x_i-x_j)^2+(y_i-y_j)^2+(z_i-z_j)^2-(\omega_i-\omega_j)^2]^{1/2}\ \
{\rm for}\ \ \kappa<0.
\end{cases}
$$
For $\kappa\ne 0$, the initial conditions must be taken such that the bodies are restricted to 3-spheres for $\kappa>0$ and hyperbolic 3-spheres for $\kappa<0$.
For $\kappa=0$, we recover the Newtonian equations,
\begin{equation}
\ddot{\bf r}_i=\sum_{j=1, j\ne i}^N\frac{m_j({\bf r}_j-{\bf r}_i)}{r_{ij}^3}, \ \ i=\overline{1,N},
\end{equation} 
with ${\bf r}_i=(x_i,y_i,z_i,0),\ i=\overline{1,N}$. The quantities $r_{ij}$ vary smoothly with $\kappa$. The values of the coordinates $\omega_i,\ i=\overline{1,N}$, and therefore the expressions $(\omega_i-\omega_j)^2, \ i,j\in\{1,2,\dots,N\}, i\ne j,$ become small when $\kappa$ gets close to 0 and vanish at $\kappa=0$. 

We further introduce the equations of motion in extrinsic coordinates and explain why they fall short of our goal (Section 2), then derive the North-Pole equations in the hope that they would solve our problem (Section 3). Unfortunately they do not, but help us get a step closer towards finding a solution. We also derive the equations of motion in intrinsic coordinates (Section 4) and explain why they also fail to address our concerns. Then we prove that all these equations can be extended to system \eqref{new}, the only framework we have found so far that offers a unified picture for all $\kappa\in\mathbb R$ (Section 5). 
We end our paper with a discussion of the bifurcations encountered by the integrals of motion when the curvature parameter $\kappa$ passes through the value $\kappa=0$  (Section 6).

\section{Equations of motion in extrinsic coordinates}
\label{extrinsic}

In this section we present the equations of motion of the curved $N$-body problem in extrinsic coordinates and explain how the flat case is obtained in the limit when $\kappa\to 0$. But before getting into details, we would like to mention why the approach of Bolyai and Lobachevsky is the natural way to extend gravitation to spaces of nonzero constant Gaussian curvature.  

The reason for introducing this extension is purely mathematical. There is no unique way of generalizing the classical equations of motion in order to recover them when the curved ambient space becomes flat. So the potential we want to use should satisfy the same basic properties the Newtonian potential does in its most basic setting---the Kepler problem---a particular case when one body moves around a fixed attracting centre. 

Two fundamental properties characterize the Newtonian potential of the Kepler problem: it is a harmonic function in 3D (but not in 2D), i.e.,\ it satisfies Laplace's equation; and it generates a central field in which all bounded orbits are closed, a result proved by Joseph Louis Bertrand in 1873, \cite{Ber}. In the early years of the 20th century, Heinrich Liebmann proved that these properties are also satisfied by the Kepler problem in spaces of constant curvature, thus offering strong arguments for 
this mathematical generalization of the gravitational force, \cite{Lieb1}, \cite{Lieb2}.

Let us further present our approach to the gravitational extension first suggested by Bolyai and Lobachevski. Take $N\ge 2$ point masses, $m_1,\dots, m_N>0$, moving on the 3-sphere (of constant Gaussian curvature $\kappa>0$),
$$
\mathbb S_\kappa^3:=\{(x,y,z,w) \ \! | \ \! x^2+y^2+z^2+w^2=\kappa^{-1}, \ \kappa>0\},
$$
viewed as embedded in $\mathbb R^4$, or on the hyperbolic 3-sphere (of constant Gaussian curvature $\kappa<0$),
$$
\mathbb H_\kappa^3:=\{(x,y,z,w) \ \! | \ \! x^2+y^2+z^2-w^2=\kappa^{-1}, \ w>0,\ \kappa<0\},
$$
viewed as embedded in the Minkowski space $\mathbb R^{3,1}$. We consider these spaces in the framework of classical mechanics, so unlike in special or general relativity, the Minkowski space mentioned above has four spatial components instead of one temporal and three spatial dimensions. So the notation $\mathbb R^{3,1}$ we adopt here rather expresses the signature of the inner product defined below instead of the nature of the components.

The coordinates of the point mass $m_i$ are given by the components of the vector ${\bf q}_i=(x_i,y_i,z_i, w_i)$, and they satisfy the constraints 
$$
x_i^2+y_i^2+z_i^2+\sigma w_i^2=\kappa^{-1},\ \ i=\overline{1,N},
$$
where $\sigma$ is the signum function
$$
\sigma:=
\begin{cases}
+1\ \ {\rm for}\ \ \kappa>0\cr
-1\ \ {\rm for}\ \ \kappa<0.
\end{cases}
$$
We define the inner product of the vectors ${\bf q}_i$ and ${\bf q}_j$ by the formula
$$
{\bf q}_i\cdot{\bf q}_j:=x_ix_j+y_iy_j+z_iz_j+\sigma z_iz_j.
$$
This is the standard inner product in $\mathbb R^4$, of signature $(+,+,+,+)$, for $\kappa>0$, but the Lorentz inner product in the Minkowski space $\mathbb R^{3,1}$, of signature $(+,+,+,-)$, for $\kappa<0$.

Let us consider the notations
$$
q^{ij}:={\bf q}_i\cdot{\bf q}_j,\  i,j\in\{1,2,\dots,N\}, \ i\ne j,
$$
$$
 q_i^2:={\bf q}_i\cdot{\bf q}_i,\ \ i=\overline{1,N}.
$$
As shown in \cite{Diacu03}, \cite{Diacu05}, or \cite{Diacu10}, the cotangent force function, 
\begin{equation}\label{force-f}
U_\kappa({\bf q})=
\begin{cases}
m_im_j\cot(d_\kappa({\bf q}_i\cdot{\bf q}_j)), \hspace{14pt} \kappa>0\cr
m_im_j\coth(d_\kappa({\bf q}_i\cdot{\bf q}_j)),\ \ \kappa<0,\cr
\end{cases}
\end{equation}
which extends the classical Newtonian force function to $\mathbb S_\kappa^3$ and $\mathbb H_\kappa^3$ for $\kappa\ne 0$ in the direction suggested by Bolyai and Lobachewski, can be put  into the form
\begin{equation}\label{forcef}
U_\kappa({\bf q})=\sum_{1\le i<j\le N}\frac{m_im_j|\kappa|^{1/2}\kappa q^{ij}}{|(\kappa q_i^2)(\kappa q_j^2)-(\kappa q^{ij})^2|^{1/2}},
\end{equation}
where
$
{\bf q}:=({\bf q}_1, {\bf q}_2,\dots, {\bf q}_N)
$
is the configuration of the particle system. But $U_\kappa$ is a homogeneous function of degree 0, so Euler's relationship,
$$
{\bf q}_i\cdot\nabla_{{\bf q}_i}U_\kappa({\bf q})=0, \ \ i=\overline{1,N},
$$
is satisfied, \cite{Euler}. Then, using the variational method of constrained Lagrangian dynamics (see \cite{Diacu03}, \cite{Diacu05}, or \cite{Diacu10}), it can be shown that the equations of motion are given by the system of differential equations
\begin{equation}\label{eqmotion}
m_i \ddot {\bf q}_i=\nabla_{{\bf q}_i} U_{\kappa}({\bf q})-\kappa m_i(\dot{\bf q}_i\cdot\dot{\bf q}_i){\bf q}_i,\ \ 
\ i=\overline{1,N},
\end{equation}
where $\kappa\ne 0$ and 
\begin{equation}\label{gradient}
{\nabla}_{{\bf q}_i}U_\kappa({\bf q})=\sum_{j=1, j\ne i}^N{m_im_j|\kappa|^{3/2}\kappa q_j^2[(\kappa q_i^2){\bf q}_j-(\kappa q^{ij}){\bf q}_i]\over
|(\kappa q_i^2)(\kappa q_j^2)-({\kappa q^{ij}
})^2|^{3/2}}, \  i=\overline{1,N}.
\end{equation}
To keep the bodies on the respective manifolds, it is enough to assume that, at the initial time $t=0$, the position vectors and the velocities satisfy the constraints
\begin{equation}\label{first-constraints}
\kappa q_i^2=1,\ \ {\bf q}_i\cdot\dot{\bf q}_i=0,
\ i=\overline{1,N},
\end{equation}
conditions that hold for all time $t$ for which the solution is defined. 

Using the constraints $\kappa q_i^2=1, \  i=\overline{1,N}$, we can write the gradient of the force function $U_\kappa$ on the manifolds of curvature $\kappa\ne 0$ as
\begin{equation}
\nabla_{{\bf q}_i}U_\kappa({\bf q})=\sum_{j=1,\ \! j\ne i}^N{{m_im_j}|\kappa|^{3/2}
\left[{\bf q}_j  - (\kappa q^{ij}){\bf q}_i \right]
\over
\left|1-\left(\kappa q^{ij}\right)^2\right|^{3/2}},  \  i=\overline{1,N}.
\label{grad}
\end{equation}

To obtain the force function \eqref{forcef}, we used in previous work the arc distance in $\mathbb S_\kappa^3$ and $\mathbb H_\kappa^3$, such that the forces between bodies acts along geodesics (see \cite{Diacu03}, \cite{Diacu05}, or \cite{Diacu07}). To include the Euclidean case, we can write in general that
$$
d_{\kappa}({\bf q}_i,{\bf q}_j):=
\begin{cases}
\kappa^{-1/2}\cos^{-1}(\kappa{\bf q}_i\cdot{\bf q}_j),\ \ \ \ \ \ \ \ \  \kappa >0\cr
|{\bf q}_i-{\bf q}_j|, \ \ \ \ \ \ \ \ \ \ \ \ \ \ \ \ \ \ \ \ \ \ \  \! \ \hspace{2 pt} \kappa=0\cr
({-\kappa})^{-1/2}\cosh^{-1}(\kappa{\bf q}_i\cdot{\bf q}_j),\hspace{10 pt} \kappa<0,\cr
\end{cases}
$$
where $d_{\kappa}({\bf q}_i,{\bf q}_j)$ is the distance between $m_i$ and $m_j$ on the manifold of curvature $\kappa\in\mathbb R$. As the curvature of the manifolds $\mathbb S_\kappa^3$ and $\mathbb H_\kappa^3$ nears 0,
the distance between $m_i$ and $m_j$ approaches the Euclidean distance, as it is obvious from geometrical considerations (see Figure \ref{common}). But this is not at all obvious from the above formula, since $|{\bf q}_i|, |{\bf q}_j|\to\infty$ as $\kappa\to 0$. So though it is intuitively clear from the geometric considerations about the distance that 
$$
\lim_{\kappa\to 0}U_\kappa({\bf q})\to U_0({\bf q}):=
\sum_{1\le i<j\le N}\frac{m_im_j}{|{\bf q}_i-{\bf q}_j|},
$$
i.e., $U_\kappa$ tends to the Newtonian force function when $\kappa\to 0$, this fact becomes less obvious when trying to use \eqref{forcef}. A similar problem appears when attempting to prove that 
$$
\lim_{\kappa\to 0}\nabla_{{\bf q}_i} U_\kappa({\bf q})=
\sum_{j=1, j\ne i}^N\frac{m_im_j({\bf q}_j-{\bf q}_i)}{|{\bf q}_j-{\bf q}_i|^3},
$$
i.e., the equations of the curved problem ($\kappa\ne 0$)
tend to the Newtonian equations when $\kappa\to 0$. But again, the above geometric considerations about the distance support the validity of this limit.

A similar technical difficulty shows up when substituting $\kappa=0$ into \eqref{eqmotion}, an operation that leads to undetermined expressions on the right hand side
of the equations of motion. So though from the geometrical and dynamical point of view we can conclude that the equations of the curved problem tend in the limit to the Newtonian equations, system \eqref{eqmotion} does not include both the curved and flat cases since the lengths of the position vectors tend to infinity when $\kappa\to 0$.
 
It is natural to suspect that the reason for this failure stays with the fact that the origin of the co-ordinate system is at the centre of the spheres and the radii of the spheres become infinite as $\kappa\to 0$. But we could shift the origin of the coordinate system to the North Pole of the 3-spheres, namely the point $(0,0,0, |\kappa|^{-1/2})$, a move that would keep the values of the coordinates finite when $\kappa\to 0$. But as we will see in the next section, this  approach alone does not fare better either.


\section{The North-Pole equations}

\begin{figure}[htbp] 
   \centering
   \includegraphics[width=2in]{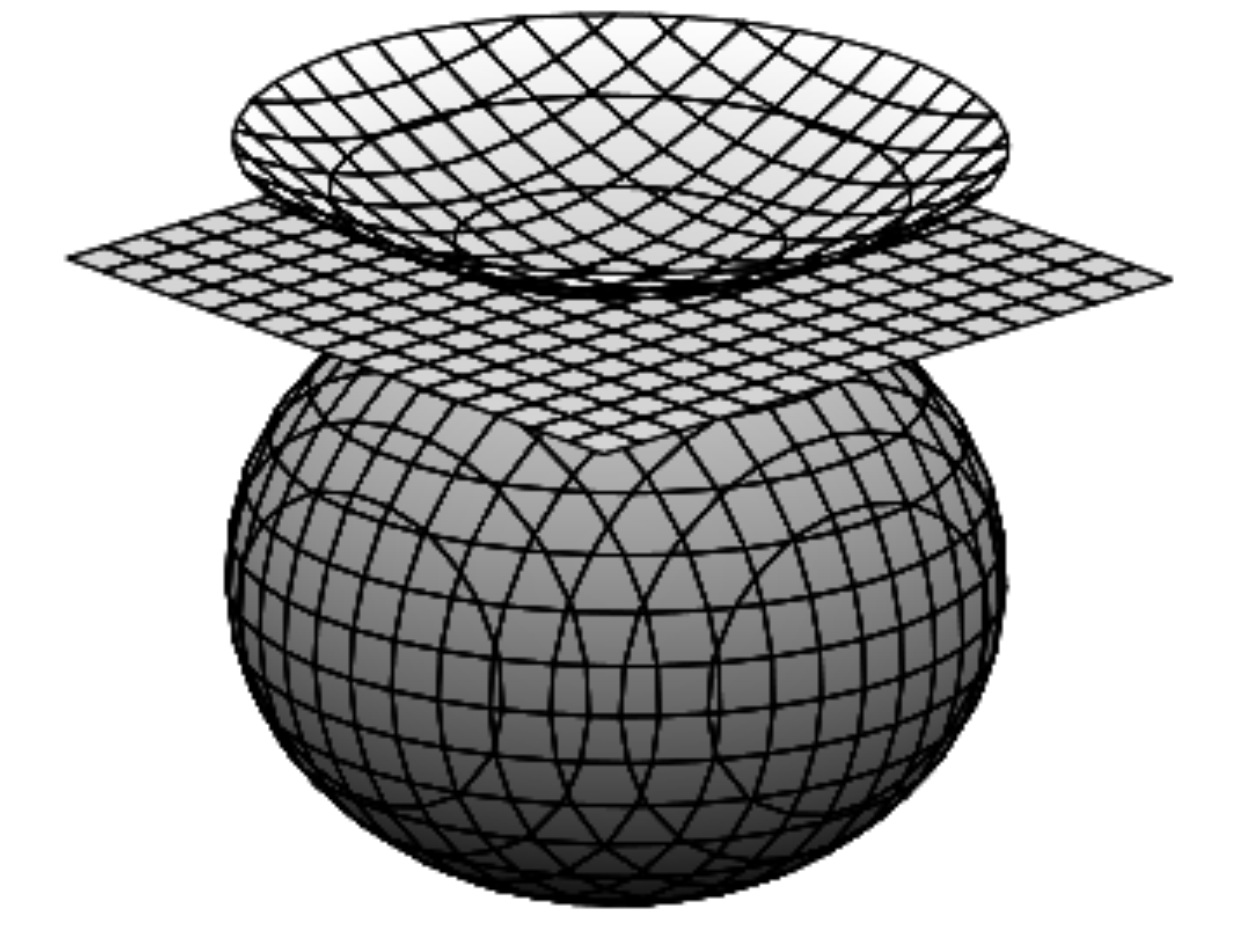}
   \caption{\small A 2D representation of the continuous transition from $\mathbb S_\kappa^3$, up, and from $\mathbb H_\kappa^3$, down, to $\mathbb R^3$. The only common point of these manifolds is the North-Pole, $(0,0,0,|\kappa|^{-1/2})$, of the 3-spheres.}
   \label{common}
\end{figure}

In this section we attempt to include the case $\kappa=0$ in the equations of motion by shifting the origin of the coordinate system to the North-Pole of the 3-spheres (see Figure \ref{common}). For this purpose we consider the change of variables
$$
\omega_i=w_i-|\kappa|^{-1/2},\ \ i=\overline{1,N},
$$
which leaves the coordinates $x_i, y_i, z_i, \ i=\overline{1,N},$ unchanged. If 
$$
\bar q^{ij}:=x_ix_j+y_iy_j+z_iz_j+\sigma\omega_i\omega_j,
$$ 
we have that
$$
\kappa q^{ij}=\kappa \bar q^{ij}+|\kappa|^{1/2}(\omega_i+\omega_j)+1,\ \ i,j\in\{1,2,\dots, N\},\ \ i\ne j.
$$

Then the equations of motion \eqref{eqmotion} take the form
\begin{equation}
\begin{cases}
\ddot x_i=\sum_{1\le i<j\le N}\frac{m_j|\kappa|^{3/2}[x_j-(\kappa\bar q^{ij}+|\kappa|^{1/2}(\omega_i+\omega_j)+1)x_i]}{|1-[\kappa\bar q^{ij}+|\kappa|^{1/2}(\omega_i+\omega_j)+1]^2|^{3/2}}-\kappa({\dot{\bar{\bf q}}}_i\cdot {\dot{\bar{\bf q}}}_i)x_i\cr
\ddot y_i=\sum_{1\le i<j\le N}\frac{m_j|\kappa|^{3/2}[y_j-(\kappa\bar q^{ij}+|\kappa|^{1/2}(\omega_i+\omega_j)+1)y_i]}{|1-[\kappa\bar q^{ij}+|\kappa|^{1/2}(\omega_i+\omega_j)+1]^2|^{3/2}}-\kappa({\dot{\bar{\bf q}}}_i\cdot {\dot{\bar{\bf q}}}_i)y_i\cr
\ddot z_i=\sum_{1\le i<j\le N}\frac{m_j|\kappa|^{3/2}[z_j-(\kappa\bar q^{ij}+|\kappa|^{1/2}(\omega_i+\omega_j)+1)z_i]}{|1-[\kappa\bar q^{ij}+|\kappa|^{1/2}(\omega_i+\omega_j)+1]^2|^{3/2}}-\kappa({\dot{\bar{\bf q}}}_i\cdot {\dot{\bar{\bf q}}}_i)z_i\cr
\ddot\omega_i=\sum_{1\le i<j\le N}\frac{m_j|\kappa|^{3/2}\{\omega_j+|\kappa|^{-1/2}-[\kappa\bar q^{ij}+|\kappa|^{1/2}(\omega_i+\omega_j)+1](\omega_i+|\kappa|^{-1/2})\}}{|1-[\kappa\bar q^{ij}+|\kappa|^{1/2}(\omega_i+\omega_j)+1]^2|^{3/2}}\cr
\hspace{3cm}-\kappa({\dot{\bar{\bf q}}}_i\cdot {\dot{\bar{\bf q}}}_i)(\omega_i+|\kappa|^{-1/2}),\ i=\overline{1,N},\cr
\end{cases}
\end{equation}
where ${\dot{\bar{\bf q}}}_i=(\dot x_i, \dot y_i, \dot z_i, \dot\omega_i)$ and 
${\dot{\bar{\bf q}}}_i\cdot {\dot{\bar{\bf q}}}_i=\dot x_i^2+\dot y_i^2+\dot z_i^2+\sigma\dot\omega_i^2,\ i=\overline{1,N}$.

As in the previous section, the equations of motion are undetermined when $\kappa=0$, although we know from the above geometrical considerations that they tend to the Newtonian equations as $\kappa\to 0$. This fact suggests that the extrinsic coordinates might not be good enough for solving our problem, so let us see if the use of intrinsic coordinates allows us to include the case $\kappa=0$ into the equations of motion.

\section{Equations of motion in intrinsic coordinates}
\label{intrinsic}

In this section we introduce the equations of motion in intrinsic coordinates in a unified context. For $\kappa<0$ and $\kappa>0$, these equations were separately derived and studied in \cite{Diacu13} and \cite{Perez}, respectively. These papers, however, treat only the 2D case. The problem of obtaining the equations of the curved $N$-body problem in intrinsic coordinates on 3D manifolds of constant curvature was not yet considered. 

So we assume in this section that the bodies move on the 2-spheres $\mathbb S_\kappa^2$ or 
the hyperbolic 2-spheres $\mathbb H_\kappa^2$, which we can write together as
$$
\mathbb M_\kappa^2=\{(\mathfrak x, \mathfrak y, \mathfrak z) \ \! |\ \! \mathfrak x^2+\mathfrak y^2+\sigma \mathfrak z^2=\kappa^{-1}, \ \kappa\ne 0, \ {\rm with}\ \ \mathfrak z>0\ \ {\rm for}\ \ \kappa<0 \}.
$$
In this new setting, the force function \eqref{forcef} has the form
\begin{equation}\label{newforcef}
\mathfrak U_\kappa({\bf p})=
\sum_{1\le i<j\le N}\frac{m_im_j|\kappa|^{1/2}\kappa p^{ij}}{|(\kappa p_i^2)(\kappa p_j^2)-(\kappa p^{ij})^2|^{1/2}},
\end{equation}
where ${\bf p}_i=(\mathfrak x_i, \mathfrak y_i, \mathfrak z_i), \ i=\overline{1,N}$, 
$$
p^{ij}:={\bf p}_i\cdot{\bf p}_j=\mathfrak x_i \mathfrak x_j+\mathfrak y_i \mathfrak y_j+\sigma \mathfrak z_i \mathfrak z_j,\ \ p_i^2:={\bf p}_i\cdot{\bf p}_i=\mathfrak x_i^2+\mathfrak y_i^2+\sigma\mathfrak z_i, \ \ i,j\in\{1,2,\dots,N\},
$$
and ${\bf p}=({\bf p}_1,{\bf p}_2,\dots, {\bf p}_N)$ is the configuration of the system. Then the equations of motion are given by
\begin{equation}\label{neweqmotion}
m_i \ddot {\bf p}_i=\nabla_{{\bf p}_i} \mathfrak U_{\kappa}({\bf p})-\kappa m_i(\dot{\bf p}_i\cdot\dot{\bf p}_i){\bf p}_i,\ \ 
\ i=\overline{1,N},
\end{equation}
where $\kappa\ne 0$, 
\begin{equation}\label{newgradient}
{\nabla}_{{\bf p}_i}\mathfrak U_\kappa({\bf p})=\sum_{\substack{j=1\\ j\ne i}}^N{m_im_j|\kappa|^{3/2}\kappa p_j^2[(\kappa p_i^2){\bf p}_j-(\kappa p^{ij}){\bf p}_i]\over
|(\kappa p_i^2)(\kappa p_j^2)-({\kappa p^{ij}
})^2|^{3/2}}, \  i=\overline{1,N},
\end{equation}
and the coordinates satisfy the constraints
$$
\kappa p_i^2=1, \ \ {\bf p}_i\cdot\dot{\bf p}_i=0,\ \ i=\overline{1,N}.
$$

To obtain the equations of motion in intrinsic coordinates, we further introduce new geometric models, both for the 2-spheres and the hyperbolic 2-spheres.  For this, we use the stereographic projection, which takes the points of coordinates $(\mathfrak x,\mathfrak y,\mathfrak z)\in$ $\mathbb M_\kappa^2$ to the points of coordinates $(u,v)$ of the plane $\mathfrak z=0$ through the bijective transformation 
\begin{equation}
\footnotesize
u=\frac{\mathfrak x}{1-\sigma|\kappa|^{1/2}\mathfrak z},\ \ v=\frac{\mathfrak y}{1-\sigma|\kappa|^{1/2} \mathfrak z}.
\label{Gprojection}
\end{equation}
The inverse of the stereographic projection takes the points of coordinates $(u,v)$ of the plane $\mathfrak z=0$ to the points
$(\mathfrak x,\mathfrak y,\mathfrak z)\in$ $\mathbb M_\kappa^2$ through the formulae
\begin{equation}
\footnotesize
\mathfrak x=\frac{2u}{1+\kappa(u^2+v^2)},\ \ \mathfrak y=\frac{2v}{1+\kappa(u^2+v^2)},\ \ \mathfrak z=\frac{\kappa(u^2+v^2)-1}{|\kappa|^{3/2}(u^2+v^2)+\sigma|\kappa|^{1/2}}.
\label{inverse}
\end{equation}
From the geometric point of view, the correspondence between a point of $\mathbb M_\kappa^2$ and
a point of the plane $\mathfrak z=0$ is made via a straight line through the point $(0,0,|\kappa|^{-1/2})$ for 
$\kappa>0$ and $(0,0,-|\kappa|^{-1/2})$ for $\kappa<0$. In the former case the projection of $\mathbb S_\kappa^2$ is $\mathbb R^2$, but with a different metric than the Euclidean one. We denote this plane by $\mathbb P_\kappa^2$. In the latter case the projection of $\mathbb H_\kappa^2$ is the Poincar\'e disk, $\mathbb D_\kappa^2$, of radius $(-\kappa)^{-1/2}$, with the corresponding hyperbolic metric. Let $\mathbb B_\kappa^2$ denote either of $\mathbb P_\kappa^2$ and $\mathbb D_\kappa^2$. With this notation we say that the stereographic projection of $\mathbb M_\kappa^2$ that preserves the geometric structure is $\mathbb B_\kappa^2$.

The metric of $\mathbb B_\kappa^2$ in coordinates $(u,v)$ is given by 
$$
ds^2=\frac{4}{[1+\kappa(u^2+v^2)]^2}(du^2+dv^2).
$$
This metric can be obtained by substituting \eqref{inverse} into 
$$
ds^2=d\mathfrak x^2+d\mathfrak y^2+\sigma d\mathfrak z^2,
$$
which defines the metric in $\mathbb R^3$, for $\sigma=1$, and in the Minkowski space ${\mathbb M}^{2,1}$, for $\sigma=-1$ (see, e.g., \cite{Dub}).
In other words, we can say that the metric in $\mathbb B_\kappa^2$ is given by the
matrix $G=(g_{ij})_{i,j=1,2}$ with
$$
g_{11}=g_{22}=\frac{4}{[1+\kappa(u^2+v^2)]^2},\ \
g_{12}=g_{21}=0.
$$
The inverse of $G$ is $G^{-1}=(g^{ij})_{i,j=1,2}$ with
$$
g^{11}=g^{22}=\frac{[1+\kappa(u^2+v^2)]^2}{4},\ \
g^{12}=g^{21}=0.
$$

Assume that the stereographic projection maps the points ${\bf q}_i$ and ${\bf q}_j$ from $\mathbb M_\kappa^2$  to the points ${\bf w}_i=(u_i,v_i)$ and ${\bf w}_j=(u_j,v_j)$ of $\mathbb B_\kappa^2$, respectively. Then, using \eqref{inverse}, we obtain that
$$
{\bf q}_i\cdot{\bf q}_j=\frac{4\kappa{\bf w}_i\cdot{\bf w}_j+(\kappa|{\bf w}_i|^2-1)(\kappa|{\bf w}_j|^2-1)}{\kappa(\kappa|{\bf w}_i|^2+1)(\kappa|{\bf w}_j|^2+1)},
$$
where ${\bf w}_i\cdot{\bf w}_j=u_iu_j+v_iv_j$, so $|{\bf w}_i|^2=u_i^2+v_i^2$.

To simplify the computations, we introduce the complex coordinates $(z,\bar z)$ with the help of the transformation
$$
z=u+iv,\ \ \bar{z}=u-iv.
$$
Then the metric of $\mathbb B_\kappa^2$ is given by
$$
ds^2=\frac{4}{(1+\kappa|z|^2)^2}\, dz\, d\bar{z},
$$
where $\frac{4}{(1+\kappa|z|^2)^2}$ is the conformal factor.

Some long but straightforward computations show that, for $\kappa\ne 0$, the above changes of variables applied to the position vectors bring the force function $U_\kappa$ given by  \eqref{newforcef} to the form
\begin{equation}
W_\kappa({\bf z},{\bf \bar{\bf z}})=\sum_{1\le i<j\le N}\frac{|\kappa|^{1/2}m_im_jB_{ij}}{|A_{ij}^2-B_{ij}^2|^{1/2}},
\end{equation}
where ${\bf z}=({z}_1,{z}_2,\dots,{z}_N),\ \bar{\bf z}=
(\bar{z}_1,\bar{z}_2,\dots,\bar{z}_N),$ and $z_i$ is the
coordinate of the body of mass $m_i, \ i=\overline{1,N}$,
$$
B_{ij}:=B(z_i,z_j, \bar{z}_i,\bar{z}_j):=2\kappa^{-1}(z_i\bar{z}_j+z_j\bar{z}_i)
+(|z_i|^2-\kappa^{-1})(|z_j|^2-\kappa^{-1}),
$$
$$
A_{ij}:=A(z_i,z_j,\bar{z}_i,\bar{z}_j):=(|z_i|^2+\kappa^{-1})(|z_j|^2+\kappa^{-1}), \ \ i,j\in\{1,2,\dots,N\}, \ i\ne j.
$$
The equations of motion \eqref{neweqmotion} take the form
\begin{equation}\label{intreqmotion}
m_i\ddot{z}_i=\frac{(\kappa|z_i|^2+1)^2}{2}\frac{\partial W_\kappa}{\partial\bar{z}_i}({\bf z},{\bf\bar z})+\frac{2|\kappa| m_i\bar{z}_i\dot{z}_i^2}{\kappa|z_i|^2+1},\ \ i=\overline{1,N},
\end{equation}
where 
$$
\frac{\partial W_\kappa}{\partial\bar{z}_i}({\bf z},{\bf\bar z})=\sum_{j=1, j\ne i}^N
\frac{2m_im_jE_{ij}}{(\sigma\kappa)^{11/2}[\sigma(A_{ij}^2-B_{ij}^2)]^{3/2}},
$$
$$
E_{ij}:=E(z_i,z_j,\bar{z}_i,\bar{z}_j):=2
(\kappa|z_i|^2+1)(\kappa|z_j|^2+1)^2
(z_j-z_i)(\kappa z_i\bar{z}_j+1).
$$

For $\kappa=0$, system \eqref{intreqmotion}  is undetermined. By looking just at these equations, it is also far from obvious that the Newtonian equations are recovered when $\kappa\to 0$, but this property is satisfied because equations \eqref{intreqmotion}  and \eqref{eqmotion} are equivalent, a result proved in 
\cite{Diacu13} and \cite{Perez}.

Since the equations of motion written in intrinsic coordinates do not solve our problem either, let us move to another attempt of finding a solution. The idea is to combine the use of extrinsic coordinates given by a frame centred at the North Pole of the 3-spheres with different distances than the geodesic ones, namely the Euclidean distance for $\kappa\ge 0$ and the Minkowski distance for $
\kappa<0$.

\section{Extension to the flat case}
\label{unification}

In this section we provide a form of the equations of motion that extends from $\kappa\ne 0$ to $\kappa=0$,
thus solving the problem we posed at the beginning of this paper. Given the position vectors ${\bf q}_i=(x_i,y_i,z_i,w_i)$ for the body $m_i$ and ${\bf q}_j=(x_j,y_j,z_j,w_j)$ for the body $m_j$, $i,j\in\{1,2,\dots, N\}, i\ne j$, let us introduce the notation
$$
q_{ij}:=
\begin{cases}
[(x_i-x_j)^2+(y_i-y_j)^2+(z_i-z_j)^2+(w_i-w_j)^2]^{1/2}\ \
{\rm for}\ \ \kappa>0\cr
[(x_i-x_j)^2+(y_i-y_j)^2+(z_i-z_j)^2]^{1/2}\ \hspace{2.52cm}
{\rm for}\ \ \kappa=0\cr
[(x_i-x_j)^2+(y_i-y_j)^2+(z_i-z_j)^2-(w_i-w_j)^2]^{1/2}\ \
{\rm for}\ \ \kappa<0.
\end{cases}
$$
For $\kappa>0$, $q_{ij}$ is the Euclidean distance between $m_i$ and $m_j$ in $\mathbb R^4$, whereas for $\kappa=0$ it represents the Euclidean distance in $\mathbb R^3$, a set that can be regarded as a hyperplane in $\mathbb R^4$. For $\kappa<0$, however, $q_{ij}$ is not a distance in the usual mathematical sense of the word. Although the quantities 
$$
(x_i-x_j)^2+(y_i-y_j)^2+(z_i-z_j)^2-(w_i-w_j)^2
$$ 
are always non-negative, such that the expressions $q_{ij}$ are positive for distinct vectors ${\bf q}_i$ and ${\bf q}_j$, it is not generally true that $q_{ik}\le q_{ij}+q_{jk}$, so this ``distance'' does not satisfy the triangle inequality. Nevertheless, it is standard to call it the Minkowski distance, although the terminology is a misnomer.

Using the fact that
$$
2q^{ij}=q_i^2+q_j^2-q_{ij}^2,
$$
which follows from a straightforward computation, the force function $U_\kappa$ given by \eqref{forcef} can be written in the ambient space as
\begin{equation}\label{gen-forcef}
V_\kappa({\bf q})=\sum_{1\le i<j\le N}\frac{m_im_j(\kappa q_i^2+\kappa q_j^2-\kappa q_{ij}^2)}{[2(\kappa q_i^2+\kappa q_j^2)q_{ij}^2-\kappa(q_i^2-q_j^2)^2-\kappa q_{ij}^4]^{1/2}}.
\end{equation}
On the manifolds of constant curvature $\kappa$, the force function $V_\kappa$ becomes
\begin{equation}
V_\kappa({\bf q})=\sum_{1\le i<j\le N}\frac{m_im_j(2-\kappa q_{ij}^2)}{q_{ij}(4-\kappa q_{ij}^2)^{1/2}},
\end{equation}
which is the same as
\begin{equation}\label{euclid-forcef}
V_\kappa({\bf q})=\sum_{1\le i<j\le N}\frac{m_im_j\Big(1-\frac{\kappa q_{ij}^2}{2}\Big)}{q_{ij}\Big(1-\frac{\kappa q_{ij}^2}{4}\Big)^{1/2}}.
\end{equation}
The dependence of $V_\kappa$ on ${\bf q}$ is obvious from the definition of the extrinsic mutual distances $q_{ij}$. We prefer to emphasize the dependence on $\bf q$ instead of the dependence on $q_{ij}$ alone because the equations of motion involve $\nabla_{{\bf q}_i}V_\kappa$. But whereas the formula of $U_\kappa$ in \eqref{forcef} cannot be extended to the flat case, the right hand side of
\eqref{euclid-forcef} makes immediate sense for $\kappa=0$. 
Since $V_\kappa$ depends only on the mutual distances, which are finite, we recover for $\kappa=0$ the classical Newtonian force function of the Euclidean space,
$$
V_\kappa({\bf q})=\sum_{1\le i<j\le N}\frac{m_im_j}{q_{ij}}.
$$

Let us now see how the equations of motion \eqref{eqmotion} get transformed. Straightforward computations show that we can put them into the form
\begin{equation}
\label{euclid}
\ddot{\bf q}_i=\sum_{j=1, j\ne i}^N\frac{m_j\Big[{\bf q}_j-\Big(1-\frac{\kappa q_{ij}^2}{2}\Big){\bf q}_i\Big]}{q_{ij}^3\Big(1-\frac{\kappa q_{ij}^2}{4}\Big)^{3/2}}-\kappa(\dot {\bf q}_i\cdot\dot{\bf q}_i){\bf q}_i, \ \  i=\overline{1,N}.
\end{equation}
For  $\kappa\ne 0$, the $2N$ initial conditions at $t=0$,
$$
\kappa q_i^2=1,\ \ \kappa{\bf q}_i\cdot\dot{\bf q}_i=0,\ \ i=\overline{1,N},
$$
must be satisfied to keep the bodies on the manifolds
$\mathbb S_\kappa^3$ or $\mathbb H_\kappa^3$.

Since the origin of the coordinate system lies at the centre of the 3-spheres, when $\kappa\to 0$ we have that $|{\bf q}_i|\to\infty$. So for $\kappa=0$ the equations are still undetermined. To overcome this last difficulty we can now make use of the idea introduced in Section 3, namely shift the origin of the coordinate system to the North-Poles $(0,0,0,|\kappa|^{-1/2})$ of the 3-spheres. For this consider again the transformation
\begin{equation}\label{transf}
\omega_i=w_i-|\kappa|^{-1/2}, \ \ i=\overline{1,N},
\end{equation}
which leaves the variables $x_i, y_i, z_i, \ i=\overline{1,N},$ unchanged, and make the notations
$$
{\bf r}=(0,0,0,\sigma|\kappa|^{1/2}),\ \ {\bf r}_i=(x_i,y_i,z_i,\omega_i),
$$
$$
r_{ij}:=
\begin{cases}
[(x_i-x_j)^2+(y_i-y_j)^2+(z_i-z_j)^2+(\omega_i-\omega_j)^2]^{1/2}\ \
{\rm for}\ \ \kappa>0\cr
[(x_i-x_j)^2+(y_i-y_j)^2+(z_i-z_j)^2]^{1/2}\ \hspace{2.45cm}
{\rm for}\ \ \kappa=0\cr
[(x_i-x_j)^2+(y_i-y_j)^2+(z_i-z_j)^2-(\omega_i-\omega_j)^2]^{1/2}\ \
{\rm for}\ \ \kappa<0.
\end{cases}
$$
By noticing that $r_{ij}=q_{ij}$, the equations of motion become
\begin{equation}\label{unified}
\ddot{\bf r}_i=\sum_{j=1, j\ne i}^N\frac{m_j\Big[{\bf r}_j-\Big(1-\frac{\kappa r_{ij}^2}{2}\Big){\bf r}_i+\frac{r_{ij}^2{\bf r}}{2}\Big]}{r_{ij}^3\Big(1-\frac{\kappa r_{ij}^2}{4}\Big)^{3/2}}-(\dot{\bf r}_i\cdot\dot{\bf r}_i)(\kappa{\bf r}_i+{\bf r}), \ \  i=\overline{1,N}. 
\end{equation}
At $t=0$, the initial conditions must have the $2N$ constraints
\begin{equation}
\label{constraints}
\kappa r_i^2+2|\kappa|^{1/2}\omega_i=0,\ \
|\kappa|^{1/2}{\bf r}_i\cdot\dot{\bf r}_i+\dot\omega_i=0,\ \ i=\overline{1,N}.
\end{equation}
Due to the invariance of $\mathbb S_\kappa^3$ and $\mathbb H_\kappa^3$, these conditions are satisfied for all $t$. They are also identically satisfied for $\kappa=0$. 

On components, system \eqref{unified} can be written as
\begin{equation}
\begin{cases}
\ddot x_i=\sum_{j=1, j\ne i}^N\frac{m_j\Big[x_j-\Big(1-\frac{\kappa r_{ij}^2}{2}\Big)x_i\Big]}{r_{ij}^3\Big(1-\frac{\kappa r_{ij}^2}{4}\Big)^{3/2}}-\kappa (\dot{\bf r}_i\cdot\dot{\bf r}_i)x_i\cr
\ddot y_i=\sum_{j=1, j\ne i}^N\frac{m_j\Big[y_j-\Big(1-\frac{\kappa r_{ij}^2}{2}\Big)y_i\Big]}{r_{ij}^3\Big(1-\frac{\kappa r_{ij}^2}{4}\Big)^{3/2}}-\kappa (\dot{\bf r}_i\cdot\dot{\bf r}_i)y_i\cr
\ddot z_i=\sum_{j=1, j\ne i}^N\frac{m_j\Big[z_j-\Big(1-\frac{\kappa r_{ij}^2}{2}\Big)z_i\Big]}{r_{ij}^3\Big(1-\frac{\kappa r_{ij}^2}{4}\Big)^{3/2}}-\kappa (\dot{\bf r}_i\cdot\dot{\bf r}_i)z_i\cr
\ddot\omega_i=\sum_{j=1, j\ne i}^N\frac{m_j\Big[\omega_j-\Big(1-\frac{\kappa r_{ij}^2}{2}\Big)\omega_i+\frac{\sigma|\kappa|^{1/2}r_{ij}^2}{2}\Big]}{r_{ij}^3\Big(1-\frac{\kappa r_{ij}^2}{4}\Big)^{3/2}}- (\dot{\bf r}_i\cdot\dot{\bf r}_i)[\kappa\omega_i+\sigma|\kappa|^{1/2}],
\end{cases}
\end{equation}
$ i=\overline{1,N}$, with the $2N$ constraints
\begin{equation}\label{last-constraints}
\begin{split}
\kappa(x_i^2+y_i^2+z_i^2+\sigma\omega_i^2)+2|\kappa|^{1/2}\omega_i=0,\hspace{1.4cm}\\
\kappa(x_i\dot x_i+y_i\dot y_i+z_i\dot z_i+\sigma\omega_i\dot\omega_i)+|\kappa|^{1/2}\dot\omega_i=0, \ \ i=\overline{1,N}.
\end{split}
\end{equation}

We can now assume that, when $\kappa$ varies, only the direction, but not the length, of the position vectors changes for given point masses on $\mathbb S_\kappa^3$ or $\mathbb H_\kappa^3$. Then, for $\kappa=0$, the values of $|{\bf r}_i|$ are finite, so we recover Newton's equations in the Euclidean case:  
\begin{equation}\label{newton}
\ddot{\bf r}_i=\sum_{j=1, j\ne i}^N\frac{m_j({\bf r}_j-{\bf r}_i)}{r_{ij}^3}, \ \ i=\overline{1,N},
\end{equation} 
where, since $\kappa=0$ and $\omega_i=0,\ i=\overline{1,N}$, the position vectors, 
$$
{\bf r}_i=(x_i,y_i,z_i,0),\ \ i=\overline{1,N},
$$ 
are free of constraints. Notice that, for consistency, we consider that the motion in $\mathbb R^3$ takes place in a hyperplane of $\mathbb R^4$, i.e., in a space of curvature $\kappa=0$ with position vectors ${\bf r}_i=(x_i,y_i,z_i,0)$ and velocities $\dot{\bf r}_i=(\dot x_i,\dot y_i,\dot z_i,0)$, so the coordinates and the velocities can be assumed to have the $2N$ constraints,
 $$
 \omega_i=\dot \omega_i=0,\ \ i=\overline{1,N},
 $$
 the same number as the constraints \eqref{constraints} that occur for $\kappa\ne 0$. Consequently the dimension of the phase space of system \eqref{newton} is $6N$, a conclusion that can be drawn either because there are no constraints in $\mathbb R^3$ or since there are $8N$ coordinate and velocity components bound by $2N$ constraints in $\mathbb R^4$.

The equations of motion \eqref{unified} are apparently less natural than the other equations of motion presented in this paper because they use the Euclidean distance in $\mathbb S_\kappa^3$ and the Minkowski distance in $\mathbb H_\kappa^3$ instead of the standard intrinsic arc distance between bodies. But for any given curvature $\kappa$, the Euclidean or the Minkowski distance uniquely determine the geodesic distance, so there is no room for confusion. Moreover, system \eqref{unified} is very convenient when we regard the classical Newtonian approach as the flat case 
of the more general problem that describes the gravitational motion of point masses in spaces of constant curvature. We emphasize that, unlike system \eqref{eqmotion}, for which it is far from obvious what happens when $\kappa\to0$, system \eqref{unified} brings forth the equations of motion of the curved $N$-body problem for any $\kappa\in\mathbb R$.

System \eqref{unified} thus opens the way towards the study of the classical $N$-body problem in the larger context of spaces of constant curvature. In particular it allows us to understand the dynamical behaviour of solutions near $\kappa=0$, an important physical problem since we actually still don't know whether the physical space is flat or curved, although it is now widely agreed that, should the curvature be nonzero, its absolute value must be very small. Although this is more of a cosmological problem, which refers to very large distances and not to those traditionally encountered in celestial mechanics, it is still an interesting mathematical problem to regard the equations describing the gravitational motion of $N$ bodies from the point of view of curved space.


\section{The integrals of motion}

In this last section we complete our paper with a study of the bifurcations that occur for the integrals of motion when the curvature parameter passes through the value $\kappa=0$. The results we obtain here show that the classical case is quite special in the context of curved space. The only integral of motion that encounters no bifurcations is the integral of energy, which exists for all $\kappa\in\mathbb R$, whereas all the other integrals change in number.  

It has been known since 1887 that the equations that describe the 3D Newtonian $N$-body problem have ten linearly independent integrals of motion that are algebraic functions relative to the position vectors and momenta and that there are no other such integrals, \cite{Bruns}. There is one integral of energy, three integrals of the centre of mass, three integrals of the linear momentum, and three integrals of the total angular momentum. As we previously proved, for $\kappa\ne 0$ there is one energy integral and six integrals of the angular momentum, but no integrals of the centre of mass and linear momentum, \cite{Diacu04}, \cite{Diacu05}, \cite{Diacu11}.  We will further show how these bifurcations occur in system \eqref{unified} when the parameter $\kappa$ passes through 0.

\subsection{The integrals of the centre of mass and the linear momentum}
The typical way to obtain the integrals of the linear momentum is to sum up $m_i\ddot{\bf r}_i$ from $i=1$ to $i=N$, notice that the obtained expression is 0, and then integrate this identity. The integrals of the centre of mass follow after another integration. More precisely, we have that
\begin{equation}
\begin{split}
\sum_{i=1}^Nm_i\ddot{\bf r}_i=\sum_{i=1}^N\sum_{j=1, j\ne i}^N\frac{m_im_j\Big[{\bf r}_j-\Big(1-\frac{\kappa r_{ij}^2}{2}\Big){\bf r}_i+\frac{r_{ij}^2{\bf r}}{2}\Big]}{r_{ij}^3\Big(1-\frac{\kappa r_{ij}^2}{4}\Big)^{3/2}}\hspace{1.3cm}\\
-\sum_{i=1}^Nm_i(\dot{\bf r}_i\cdot\dot{\bf r}_i)(\kappa{\bf r}_i+{\bf r})=\hspace{3.5cm}\\
\sum_{i=1}^N\sum_{j=1, j\ne i}^N\frac{m_im_j\frac{r_{ij}^2}{2}(\kappa{\bf r}_i+{\bf r})}{r_{ij}^3\Big(1-\frac{\kappa r_{ij}^2}{4}\Big)^{3/2}}-\sum_{i=1}^Nm_i(\dot{\bf r}_i\cdot\dot{\bf r}_i)(\kappa{\bf r}_i+{\bf r}), \hspace{1.1 cm}
\end{split}
\end{equation}
which is 0 for any solution only if $\kappa=0$. By integrating in the case $\kappa=0$, we obtain the three integrals of the linear momentum,
\begin{equation}\label{linmom}
\sum_{i=1}^Nm_i\dot{\bf r}_i={\bf a},
\end{equation}
where ${\bf a}=(a_1,a_2,a_3)$ is an integration vector.
By integrating equations \eqref{linmom}, we are led to the integrals of the centre of mass,
\begin{equation}\label{centrmass}
\sum_{i=1}^Nm_i{\bf r}_i-{\bf a}t={\bf b},
\end{equation}
where ${\bf b}=(b_1,b_2,b_3)$ is another integration vector. Obviously, these integrals do not show up for $\kappa\ne 0$, a fact that puts into the evidence the bifurcations these integrals encounter at $\kappa=0$. 

From the dynamical point of view, the integrals \eqref{linmom} and \eqref{centrmass} express the fact that
the centre of mass of the particle system moves uniformly along a straight when ${\bf a}\ne {\bf 0}$. By taking the origin of the coordinate system at the centre of mass, which implies that ${\bf a}={\bf b}={\bf 0}$, the above integrals become, respectively,
\begin{equation}\label{linmom2}
\sum_{i=1}^Nm_i\dot{\bf r}_i={\bf 0},
\end{equation} 
\begin{equation}\label{centrmass2}
\sum_{i=1}^Nm_i{\bf r}_i={\bf 0}.
\end{equation}
Their physical interpretation is that the centre of mass is fixed relative to the coordinate system. This means that the
forces acting on the centre of mass cancel each other. In general, no such physical properties occur when $\kappa\ne 0$. In particular, there is no point at which the forces acting on it cancel each other. Nevertheless, some particular solutions of the equations of motion have this property,
as shown in previous work, such as \cite{Diacu03} and \cite{Diacu05}.

\subsection{The integral of energy}
We further obtain the integral of energy for system \eqref{euclid} and then use the change of variables \eqref{transf} to derive this integral for system \eqref{unified}. The standard approach is to take $m_i\ddot{\bf q}_i\cdot\dot{\bf q}_i$ and sum up from $i=1$ to $i=N$, i.e.,
$$
\sum_{i=1}^Nm_i\ddot{\bf q}_i\cdot\dot{\bf q}_i=
\sum_{i=1}^N\dot{\bf q}_i\cdot\nabla_{{\bf q}_i}V_\kappa({\bf q})-\sum_{i=1}^Nm_i(\dot{\bf q}_i\cdot\dot{\bf q}_i)(\kappa{\bf q}_i\cdot\dot{\bf q}_i)=\frac{d}{dt}V_\kappa({\bf q}).
$$
By integration we obtain the energy integral, 
\begin{equation}
H_\kappa({\bf q}, \dot{\bf q}):=T_\kappa({\bf q}, \dot{\bf q})-V_\kappa({\bf q})=h,
\end{equation}
where $H_\kappa$ is the Hamiltonian function, 
$$
T_\kappa({\bf q}, \dot{\bf q}):=\frac{1}{2}\sum_{i=1}^N\kappa m_iq_i^2(\dot{\bf q}_i\cdot\dot{\bf q}_i)
$$
is the kinetic energy, and $h$ is an integration constant. Using the transformations \eqref{transf}, the kinetic energy $T_\kappa$ becomes 
$$
\mathcal T_\kappa({\bf r},\dot{\bf r})=
\frac{1}{2}\sum_{i=1}^Nm_i(\kappa r_i^2+2|\kappa|^{1/2}\omega_i+1)(\dot{\bf r}_i\cdot\dot{\bf r}_i),
$$
so the integral of energy for system \eqref{unified} takes the form
\begin{equation}
\frac{1}{2}\sum_{i=1}^Nm_i(\kappa r_i^2+2|\kappa|^{1/2}\omega_i+1)(\dot{\bf r}_i\cdot\dot{\bf r}_i)
-\sum_{1\le i<j\le N}\frac{m_im_j\Big(1-\frac{\kappa r_{ij}^2}{2}\Big)}{r_{ij}\Big(1-\frac{\kappa r_{ij}^2}{4}\Big)^{1/2}}=h.
\end{equation}
For $\kappa=0$, we recover the well-known integral of the Newtonian equations,
$$
\frac{1}{2}\sum_{i=1}^Nm_i({\dot x}_i^2+{\dot y}_i^2+{\dot z}_i^2)-\sum_{1\le i<j\le N}\frac{m_im_j}{r_{ij}}=h,
$$
so no bifurcations occur in this case.

\subsection{The integrals of the total angular momentum}
As in the case of the energy integral, we derive the integrals of the total angular momentum for system \eqref{euclid} and use the change of variables \eqref{transf} to obtain the integrals for equations  \eqref{unified}. The total angular momentum is defined as
$$
\sum_{i=1}^Nm_i{\bf q}_i\wedge\dot
{\bf q}_i,
$$ 
where $\wedge$ represents the exterior product of the
Grassman algebra over ${\mathbb R}^4$. We show that this quantity is conserved for the equations of motion \eqref{euclid}, i.e.,
\begin{equation}
\sum_{i=1}^Nm_i{\bf q}_i\wedge\dot{\bf q}_i={\bf c},
\label{angintegrals}
\end{equation}
where
${\bf c}=c_{wx}{\bf e}_w\wedge {\bf e}_x+c_{wy}{\bf e}_w\wedge {\bf e}_y+c_{wz}{\bf e}_w\wedge {\bf e}_z+
c_{xy}{\bf e}_x\wedge {\bf e}_y+c_{xz}{\bf e}_x\wedge {\bf e}_z+c_{yz}{\bf e}_y\wedge {\bf e}_z,$ with the coefficients 
$c_{wx}, c_{wy}, c_{wz}, c_{xy}, c_{xz}, c_{yz}\in{\mathbb R}$, and
$$
{\bf e}_x=(1,0,0,0),\ {\bf e}_y=(0,1,0,0),\ {\bf e}_z=(0,0,1,0),\ {\bf e}_w=(0,0,0,1)
$$
representing the vectors of the standard basis of $\mathbb R^4$. We obtain this conservation law by integrating the identity formed by the left and right expressions in the sequence of equations
\begin{equation}
\begin{split}
\sum_{i=1}^Nm_i\ddot{\bf q}_i\wedge{\bf q}_i=\sum_{i=1}^N\sum_{j=1, j\ne i}^N\frac{m_im_j{\bf q}_j\wedge{\bf q}_i}{q_{ij}^3\Big(1-\frac{\kappa q_{ij}^2}{4}\Big)^{3/2}}
\hspace{1.8cm}\\
-\sum_{i=1}^N\left[\frac{m_im_j(1-\frac{\kappa q_{ij}^2}{2}\Big)}{q_{ij}^3\Big(1-\frac{\kappa q_{ij}^2}{4}\Big)^{3/2}}-
\kappa m_i(\dot {\bf q}_i\cdot\dot{\bf q}_i) \right]{\bf q}_i\wedge{\bf q}_i
={\bf 0},\hspace{1.2cm}
\end{split}
\end{equation}
which follows after $\wedge$-multiplying the equations of motion \eqref{euclid} 
by $m_i{\bf q}_i$ and summing up from $i=1$ to $i=N$. The last of the above identities follows from the skew-symmetry of the $\wedge$ operation and, consequently, from the fact that  ${\bf q}_i\wedge{\bf q}_i={\bf 0}, \ i=\overline{1,N}$.
On components, the six integrals in \eqref{angintegrals} can be written as

\begin{equation}
\sum_{i=1}^Nm_i(x_i\dot{y}_i-\dot{x}_iy_i)=c_{xy},
\end{equation}
\begin{equation}
\sum_{i=1}^Nm_i(x_i\dot{z}_i-\dot{x}_iz_i)=c_{xz},
\end{equation}
\begin{equation}
\sum_{i=1}^Nm_i(y_i\dot{z}_i-\dot{y}_iz_i)=c_{yz},
\end{equation}
\begin{equation} 
\sum_{i=1}^Nm_i(w_i\dot{x}_i-\dot{w}_ix_i)=c_{wx},
\end{equation}
\begin{equation} 
\sum_{i=1}^Nm_i(w_i\dot{y}_i-\dot{w}_iy_i)=c_{wy},
\end{equation}
\begin{equation}
\sum_{i=1}^Nm_i(w_i\dot{z}_i-\dot{w}_iz_i)=c_{wz}.
\end{equation} 
Using the transformations \eqref{transf}, we can see that for system \eqref{unified} these integrals take the form
\begin{equation}
\sum_{i=1}^Nm_i(x_i\dot{y}_i-\dot{x}_iy_i)=c_{xy},
\end{equation}
\begin{equation}
\sum_{i=1}^Nm_i(x_i\dot{z}_i-\dot{x}_iz_i)=c_{xz},
\end{equation}
\begin{equation}
\sum_{i=1}^Nm_i(y_i\dot{z}_i-\dot{y}_iz_i)=c_{yz},
\end{equation}
\begin{equation}
\sum_{i=1}^Nm_i{\dot x}_i+|\kappa|^{1/2}\sum_{i=1}^Nm_i(\omega_i\dot{x}_i-\dot{\omega}_ix_i)=|\kappa|^{1/2}c_{wx},
\end{equation} 
\begin{equation}
\sum_{i=1}^Nm_i{\dot y}_i+|\kappa|^{1/2}\sum_{i=1}^Nm_i(\omega_i\dot{y}_i-\dot{\omega}_iy_i)=|\kappa|^{1/2}c_{wy},
\end{equation}
\begin{equation}
\sum_{i=1}^Nm_i{\dot z}_i+|\kappa|^{1/2}\sum_{i=1}^Nm_i(\omega_i\dot{z}_i-\dot{\omega}_iz_i)=|\kappa|^{1/2}c_{wz}. 
\end{equation}
So for $\kappa\ne 0$, system \eqref{unified} has six integrals of the total angular momentum. But at $\kappa=0$, only the first three integrals of the total angular momentum remain; the others become the three integrals of the linear momentum  \eqref{linmom2} obtained when the origin of the coordinate system is taken at the centre of mass of the particle system. So an interesting kind of bifurcation occurs in this case as we pass through the value $\kappa=0$ of the curvature parameter. 


\section{Conclusions}

The study of the new equations of motion \eqref{unified} is neither simpler nor more complicated than that of the systems provided in intrinsic or extrinsic coordinates in Sections 2, 3, and 4, although certain problems might be more approachable in the framework of some equations than in the context of others. Nevertheless, system \eqref{unified} has the advantage of unifying the cases $\kappa\ne 0$ and $\kappa=0$, thus offering a larger perspective for the Newtonian equations of the $N$-body problem. As we showed above, some interesting bifurcations occur for the integrals of motion as we pass through the value $\kappa=0$ of the curvature parameter.

All the solutions obtained for $\kappa=0$ can now be viewed from the point of view of $\kappa\ne 0$ and vice versa, in the sense of studying whether such solutions bifurcate or occur for all values of the curvature parameter $\kappa$, under what circumstances they show up, and whether the stability of these solutions changes with $\kappa$. Consequently system \eqref{unified} opens new perspectives of research that were not possible with the previously derived equations of motion. In particular, the study of the Lagrangian and the Eulerian orbits of the 3-body problem can be considered in the future in the context of this larger framework. But these are a couple among the many exciting problems that can be regarded from this novel point of view.


\end{document}